\documentclass{amsart}
\usepackage{amsfonts,amssymb,amsmath,amsthm}
\usepackage{url}
\usepackage{enumerate}

\newtheorem{thrm}{Theorem}[section]
\newtheorem{lem}[thrm]{Lemma}
\newtheorem{prop}[thrm]{Proposition}
\newtheorem{cor}[thrm]{Corollary}
\theoremstyle{definition}
\newtheorem{definition}[thrm]{Definition}
\newtheorem{remark}[thrm]{Remark}

\newtheorem{question}[thrm]{Question}
\newtheorem{Note}[thrm]{Notification}
\numberwithin{equation}{section}
\author{M. Azizi and A. Rejali}
\address{Department of Mathematics\\
University of Isfahan\\
Isfahan, Iran} \email{m.azizi@sci.ui.ac.ir}
\address{Department of Mathematics\\
University of Isfahan\\
Isfahan, Iran} \email{rejali@sci.ui.ac.ir}



\begin{document}
\title[Urysohn and Tietze extensions of Lipschitz functions]
{Urysohn and Tietze extensions of Lipschitz functions} \maketitle

\begin{abstract}
Let $(X,d)$ be a metric space and $\alpha>0$. In this paper, we
study extensions of some complex-valued Lipschitz functions, from
some special subset $X_0$ to $X$. These extensions are with
no-increasing Lipschitz number or the smallest Lipschitz number.
Moreover, we show that under some conditions, Tietze extension
theorem can be generalized for Lipschitz functions and call it
Tietze-Lipschitz extension. Furthermore, we generalize
Urysohn-lemma for Lipschitz functions. In fact we present a
necessary and sufficient condition for that Lipschitz functions
separate subsets of $X$.
\end{abstract}

\section{Introduction and Preliminaries}

Let $(X,d)$ be a metric space and $\alpha>0$. Then $Lip_\alpha
{X}$ is the subspace of $B(X)$, consisting of all bounded
complex-valued functions $f$ on $X$, such that
\begin{equation}\label{e13}
p_{\alpha}(f):=\sup\left\{\frac{|f(x)-f(y)|}{d(x,y)^{\alpha}}:
x,y\in X,\;\;x\neq y\right\}<\infty.
\end{equation}
It is known that $Lip_\alpha {X}$, endowed with the norm
$\|.\|_{\alpha}$ given by
$$
\|f\|_\alpha=p_\alpha(f)+\|f\|_{\infty},
$$
is a Banach space, called Lipschitz space. We call
$p_{\alpha}(f)$, Lipschitz number of $f$. As the notation, used in
\cite{W}, we use $L(f)$ instead of $p_{1}(f)$. Note that all
Lipschitz functions are uniformly continuous but the converse of
this statement is not necessarily true; see \cite[Page 3]{W}. It
is obvious that all Lipschitz spaces contain $Cons(X)$, the space
containing of all complex-valued constant function on $X$.

Similar to continuous functions, combining two Lipschitz functions
(if possible) is again a Lipschitz function; see \cite[Page 5]{W}
for complete information regarding Lipschitz functions and more
generally the Lipschitz morphism between two metric spaces.

Lipschitz spaces were first considered by Sherbert \cite{S}; see
also Bishop \cite{B}. Also Lipschitz spaces has been considered by
some recent authors; for example \cite{aalr4}, \cite{DNS},
\cite{FN}, \cite{HMT} and \cite{MW}. In \cite{aalr4}, we did a
detailed study, concerning the structure of Lipschitz spaces.
Moreover, we investigated arbitrary intersections of Lipschitz
spaces, denoted by $\cap_{\alpha\in J}Lip_{\alpha}X$, where $J$ is
an arbitrary subset of $(0,\infty)$. Furthermore, we considered
Lipschitz spaces as Banach algebras, associated with pointwise
product and studied $C-$character amenability of Lipschitz algebra
$Lip_\alpha X$.

In \cite{W}, it has been studied another important issue related
to Lipschitz functions. In fact it has been investigated when a
real or complex-valued Lipschitz function has an extension from
the subset $X_0\subseteq X$ to $X$, without increasing Lipschitz
number, or with the least changes in Lipschitz number. Most of the
results are for the case where $\alpha=1$. In most situations
there is not important differences between $\mathbb{R}$ and
$\mathbb{C}$. But there exists some situations that complex
version is harder than real version. Namely by \cite[Theorem
1.5.6]{W}, every real-valued function $f_0\in Lip_1X_0$ has an
extension $f$ to $X$ such with $L(f)=L(f_0)$ and
$\|f_0\|_\infty=\|f\|_\infty$. But by part (b) of \cite[Theorem
1.5.6]{W}, if $f_0$ is complex-valued then there exists
$f:X\rightarrow\Bbb C$ such that $L(f)\leq\sqrt{2}\;L(f_0)$. In
fact the extension $f$ with $L(f)=L(f_0)$ cannot be necessarily
achieved in the complex case; see \cite[Example 1.5.7]{W}. In this
example, $X_0$ is a three point set and the range of $f_0$ is an
equilateral triangle . Also $X=X_0\cup\{e\}$. It has been defined
an extension $f$ of $f_0$ to $X$, with the smallest Lipschitz
number $\dfrac{2}{\sqrt{3}}$. This example shows that in general
one cannot achieve extensions with the Lipschitz number less than
$\dfrac{2}{\sqrt{3}}$. Note that if $f_0:X_0\to \mathbb{C}$ is a
Lipschitz function, $e\in X-X_0$ and if $f$ is an extension of
$f_0$ to $X_0\cup\{e\}$, then we have always
$\dfrac{L(f)}{L(f_0)}\geq 1$. It is noticeable that this fraction
is nearer to $1$. In particular extensions with constant $1$, are
very well. As it is mentioned in \cite[Page 19]{W}, this is known
from Kirszbraun`s theorem \cite{F} that if $X\subset {\mathbb
R}^n$ is equipped with the inherited Euclidean metric, $X_0\subset
X$ and $f_0:X_0\to \mathbb{C}$ is a Lipschitz function, then it is
possible to extend $f_0$ to $X$, without increasing its Lipschitz
number. Mainly, by \cite[Page 19]{W}, complex-valued Lipschitz
functions on subsets of $X$ extend to whole of $X$ without
increasing Lipschitz number if and only if every four-element
subset of $X$ has this property.

After all this discussions, an open problem has been formulated in
\cite[page 18]{W} as the following:

\begin{question}\label{q1}
What is the best possible constant in part (b) of \cite[Theorem
1.5.6]{W}?
\end{question}

The second section of the present work is devoted to rather
investigate Question \ref{q1}. In fact suppose that $X$ is a set
and $X_0=X-\{e\}$, and $d$ is a metric on $X$ defined as:
\begin{eqnarray*}
 d(x,y)= \left\lbrace
 \begin{array}{c l}
  0& \hbox{\;\;  if $x=y$}\\
  1& \hbox{\;\;  if $x\neq y$}.
 \end{array}\right.
\end{eqnarray*}
and $d(e,x)=l$ for all $x,y\in X_0$, where $l$ is a fixed number
in $[\dfrac{1}{2},\infty)$. As the main result we prove that every
triangle $f_0$ on $X_0=X-\{e\}$ has an extension $f$ to $X$ such
that
$$
L(f)\leq\frac{2}{\sqrt{3}}L(f_0).
$$
We also show that $\dfrac{2}{\sqrt{3}}$ is the best (least)
increasing coefficient. We also investigate the conditions such
that there exists an extension with no increasing in Lipschitz
number. This is in fact a generalization of part (a) of
\cite[Theorem 1.5.6]{W}, to complex version. Indeed, as a main
result we show that $f_0$ has an angle that is greater than or
equal to $\sin^{-1}(\dfrac{1}{2l})$ if and only if $f_0$ has an
extension $f$ on $X$  such that $L(f)=L(f_0)$. Moreover, we extend
our main results to some special tetragon and also regular
$n-$polygons.

In section 3, we generalize Tietze extension theorem for
complex-valued Lipschitz functions. In fact we show that under
some conditions every $f\in Lip_\alpha X_0$ $(X_0\subseteq X)$,
can be extended to a function $f\in Lip_\alpha X$, preserving
Lipschitz norm.

The last section contains the Lipschitz version of Urysohn`s
lemma. Let $(X,d)$ be a metric space and $0<\alpha\leq 1$. In
\cite[page 4]{W}, it has been claimed that for all disjoint closed
subsets $A$ and $B$ of $X$, there is a Lipschitz function
$f:X\rightarrow\Bbb R$, constantly zero on $A$ and constantly 1 on
$B$. It has been mentioned that the existence of this function is
a kind of metric version of Urysohn's Lemma. In fact it has been
assumed that $d(A,B)>0$, for all disjoint closed subsets $A$ and
$B$ of $X$. But there is a gap in this result. We provide an
example to show that this statement is not necessarily true.
Finally, as the main result we prove that if $A$ and $B$ are two
subsets of the metric space $(X,d)$ then $Lip_{\alpha}X$ separates
them if and only if $d(A,B)>0$.

\section{\bf Extensions of Lipschitz functions with the best Lipschitz number}

In this section, we present the best Lipschitz number of some
Lipschitz functions, for some special metric spaces. We first
introduce essential concepts and terms, which will be require
throughout the paper.

\begin{definition}
Let $X$ be a nonempty set.
\begin{enumerate}
\item Suppose that $f:X\to \mathbb{C}$ is a $n-$valued function on
$X$ with range $\{p_1,\cdots,p_n\}$. Then $f$ is called an
"$n$-polygon" on $X$ if $\{p_1,\cdots,p_n\}$ are the vertices of
an $n$-polygon in $\mathbb{C}$. In particular, $f$ is called a
"triangle", if $n=3$. \item Let $f,g$ be two triangles on $X$. We
say that $f$ and $g$ are "congruent", if their ranges are
congruent as two triangles with respect to three sides. In this
case we write $f\approx g$.
\end{enumerate}
\end{definition}

Throughout the paper, the largest side of the triangle $f$ is
called the base of the triangle. Also the middle point of the base
of $f$ is denoted by $M_f$. Moreover we denote by $C_f$, the
common point that all perpendicular of all sides of a $n-$polygon,
cross each other (if available) and we call it the center of $f$.
Moreover we denote by $E$, the distance between $C_f$ and each
vertex of the triangle. Note that every regular $n-$polygon can be
inscribed in a circle and $C_f$ is exactly the center of the
circle. Furthermore the distance between $C_f$ and each of the
vertices is a unique number. Note that every triangle (not
necessarily regular) has the unique center.

\begin{definition}
Let $f$ be a triangle. Then $f$ is called normal if its range has
the following properties:
\begin{enumerate}
\item all vertices are in the first quarter of $\mathbb{C}$, \item
one vertex is the element zero, \item the largest side of $f$ lies
on the $x-$axis, \item one vertex that is located in front of the
base, is not placed in the right side of the perpendicular of the
base. This point is called index point of $f$ and is denoted by
$D_f$.
\end{enumerate}
\end{definition}

Let $f$ be a triangle on the set $X$, with the vertices
$\{p_1,p_2,p_3\}$. Also suppose that $A_i$ is the inverse image of
$\{p_i\}$ by $f$. Using rotation, translation and reflection
actions, one can find a normal triangle which is congruent to the
range of $f$. We show its vertices by $\{q_1,q_2,q_3\}$, where
$q_i$ is the point which is obtained by actions on $p_i$. Now
define the function $f^\sharp$ on $X$ as $f^\sharp(A_i)=q_i$. It
is easily verified that $f\approx f^\sharp$. We say to
$f^{\sharp}$, normalized of $f$. Thus the following essential
proposition is obtained, immediately.

\begin{prop}\label{lem}
Let $(X,d)$ be a metric space, $f$ be a triangle on $X$ and
$\alpha>0$. Then $p_\alpha(f)=p_\alpha(f^\sharp)$.
\end{prop}

The following results are very useful in our main result. Both are
known in elementary geometry.

\begin{prop}\label{right angle}
Let $(X,d)$ be a metric space and $f$ be a normal triangle on $X$.
Suppose that the largest angle and the length of base of $f$ are
$\theta$ and $r$, respectively. Then $\theta\leq \dfrac{\pi}{2}$
if and only if $C_{f}$ is above or over the base of $f$.
\end{prop}

\begin{proof}
Suppose that $A=0$, $B=r$ and $C=x+iy$ are the vertices of $f$.
First suppose that $\theta\leq\dfrac{\pi}{2}$. Then $C$ is out of
the circle, with the center $\dfrac{r}{2}$ and radius
$\dfrac{r}{2}$. Thus
$$\frac{r}{2}\leq |C-\frac{r}{2}| =\sqrt{(x-\frac{r}{2})^2+y^2}.$$
It follows that $x^2-rx+y^2\geq 0$. Also since $f$ is normal then
$y>0$ and so $y+\dfrac{x(x-r)}{y}\geq 0$. Consequently
$C_{f}=\dfrac{1}{2}(r+i(\frac{x(x-r)}{y}+y))$ is located above or
over the $x$-axis, as claimed. The reverse of the implication is
obvious.
\end{proof}

\begin{prop}\label{distance}
Let $(X,d)$ be a metric space and $f$ be a normal triangle on $X$.
Suppose that the largest angle and the length of base of $f$ are
$\theta$ and $r$, respectively. Then
$$
E=\frac{r}{2\sin\theta}.
$$
\end{prop}

\begin{proof}
Suppose that the vertices of triangle are, $A=0$, $B=r$ and $C$
such that $\hat{C}=\theta$, where $\hat{C}$ is the angle of the
vertex $C$. Also take $O:=C_{f}$, the center of the triangle and
so the center of peripheral circle of $f$. Now Consider two
following cases:
\begin{enumerate}
\item If $\theta\leq \dfrac{\pi}{2}$ then by Lemma \ref{right
angle}, the center of the triangle is located above or over the
base of $f$. It follows that $\hat{AOB}=2\hat{C}$. But triangle
$AOB$ is obviously isosceles. If we denote perpendicular leg of
the side $AB$ by $H$, then we have $\hat{AOH}=\hat{C}=\theta$ and
so
$$
\sin\theta=\frac{AH}{AO}=\frac{r/2}{AO}.
$$
It follows that
$$
E=\frac{r}{2\sin\theta}.
$$
\item If $\theta> \frac{\pi}{2}$, then by Lemma \ref{right angle},
the center of the triangle is under the base of $f$. Some
arguments analogous to those given in part (1), we obtain
$$
\sin(\pi-\theta)=\frac{AH}{AO}=\frac{r/2}{AO},
$$
and so again $E=\dfrac{r}{2\sin\theta}$.
\end{enumerate}
\end{proof}

We also require the following essential lemma. The proof is
straightforward and so is left to the reader.

\begin{lem}\label{3}
Let $l_1$ and $l_2$ be two lines in the complex plane,
intersecting at $M$. Moreover suppose that $A$ is one of the four
parts, made by these lines, such that $p\notin A$.
\begin{enumerate}
\item[(i)] If it can be straightened from $p$ to $A$, and the
perpendicular leg is $H$, then
$$ |p-A|=\inf\{|p-x|:x\in A\}=|p-H|.
$$
\item[(ii)] If it cannot be straightened from $p$ to $A$, then
$$
|p-A|=\inf\{|p-x|:x\in A\}=|p-M|.
$$
\end{enumerate}
\end{lem}

We are now in a position to prove the main result of this section.

\begin{thrm}\label{4}
Let $X$ be a set with at least three elements, $e\in X$ and
$X_0=X-\{e\}$. Suppose that $d$ is a metric on $X$, defined as:
\begin{eqnarray*}
 d(x,y)= \left\lbrace
 \begin{array}{c l}
  0& \hbox{\;\;  if $x=y$}\\
  1& \hbox{\;\;  if $x\neq y$}.
 \end{array}\right.
\end{eqnarray*}
and $d(e,x)=l$ for all $x,y\in X_0$, where $l$ is a fixed number
in $[\dfrac{1}{2},\infty)$. Moreover suppose that $f_0$ is a
triangle on $X_0$. Then $f_0$ has an extension $f$ to $X$ such
that
\begin{equation}\label{2}
L(f)\leq \frac{2}{\sqrt{3}}L(f_0).
\end{equation}
Moreover the best (least) possible constant for the above
inequality is $\dfrac{2}{\sqrt{3}}$.
\end{thrm}

\begin{proof}
We should only define $f(e)$. Suppose that the vertices of
${f_0}^\sharp$ are $A:=0$ ,$B:=r$ and $C$. Moreover suppose that
$r$ is the length of the base of ${f_0}^\sharp$. Thus
$L({f_0}^\sharp)=r$, obviously. Also $C$ can be changed in the set
$D$, defined as the following:
$$D:=\{x+iy\in \mathbb{C}:|x+iy|\leq r , |(x+iy)-r|\leq r , y>0 , x\leq r/2\}.$$
Note that $D$ is the union of two following disjoint sets:
$$D_1:=\{z\in D:|z-\frac{r}{2}|\leq \frac{r}{2}\},$$
and
$$D_2:=\{z\in D:|z-\frac{r}{2}|> \frac{r}{2}\}.$$
Now define
\begin{eqnarray*}
  f_{x,y}(e):= \left\lbrace
  \begin{array}{c l}
  M_{{f_0}^\sharp}= \frac{1}{2}r& \hbox {\;\; if $C=x+iy \in D_1$}\\
  C_{{f_0}^\sharp}= \frac{1}{2}(r+i(\frac{x(x-r)}{y}+y))& \hbox{\;\;  if $C=x+iy\in D_2$}.
  \end{array}\right.
\end{eqnarray*}
If $C\in D_1$ then we have
$$
L(f_{x,y})=\max\left\{L({f_0}^\sharp),\max\{\frac{r}{2l},\frac{1}{l}|\frac{r}{2}-C|\}\right\}.
$$
Since $C\in D_1$, thus $|r/2-C|\leq r/2$ and so
\begin{equation}\label{5}
L(f_{x,y})=\max\{r,r/2l\}=r=L({f_0}^\sharp).
\end{equation}
Now suppose that $C\in D_2$. In this case,
$\hat{C}=\theta<\frac{\pi}{2}$. Since the distance between
$C_{{f_0}^{\sharp}}$ and all vertices are the same, it follows
from lemma \ref{distance} that
\begin{eqnarray*}
L(f_{x,y})&=&\max\left\{L(f_0),\max\{\frac{1}{l}|f_{x,y}(e)|\;:\;C=x+iy\in
D_2\}\right\}\\&=&\max\left\{r,\max\{\frac{r}{2l\sin\theta}\;:\;\frac{\pi}{3}\leq
\theta< \frac{\pi}{2}\}\right\}\\&=&\left\lbrace
 \begin{array}{c l}
  r& \hbox{\;\;  if $l\geq \frac{\sqrt{3}}{3}$}\\
  \frac{r}{\sqrt{3}l}& \hbox{\;\;  if $l< \frac{\sqrt{3}}{3}$}
 \end{array}\right.\\&=&\max\{r,\frac{r}{\sqrt{3}l}\}\\ &\leq &\frac{2}{\sqrt{3}}r\\&=&\frac{2}{\sqrt{3}}L({f_0}^{\sharp}).
\end{eqnarray*}
Note that the function $\dfrac{r}{2l\sin\theta}$ is decreasing
with respect to $\theta$. Moreover since $\theta$ is the greatest
angle in the triangle, it follows that $\theta\geq\dfrac{\pi}{3}$.
Thus we obtain
\begin{equation}\label{6}
L(f_{x,y})\leq \frac{2}{\sqrt{3}}L({f_0}^{\sharp}).
\end{equation}
The implications \eqref{5} and \eqref{6} imply that $f:=f_{x,y}$
is the desired extension of ${f_0}^{\sharp}$. Now we prove the
statement for $f_0$. Define
\begin{eqnarray*}
  f(e):= \left\lbrace
  \begin{array}{c l}
  M_{f_0} & \hbox {\;\; if ${f}^{\sharp}(e)=M_{{f_0}^{\sharp}}$ }\\
  C_{f_0} &  \hbox{\;\; if ${f}^{\sharp}(e)=C_{{f_0}^{\sharp}}$},
  \end{array}\right.
\end{eqnarray*}
where ${f}^{\sharp}$ is the extension of ${f_0}^{\sharp}$, which
there exists by the above discussion. By Proposition \ref{lem}, we
obtain
$$
L(f)=L({f}^{\sharp})\leq
\frac{2}{\sqrt{3}}L({f_0}^{\sharp})=\frac{2}{\sqrt{3}}L(f_0).
$$
Thus $f$ is the desired extension of $f_0$. Now we show that the
best (least) possible constant for the inequality \eqref{2} is
$\dfrac{2}{\sqrt{3}}$. Suppose that $l=1/2$ and $f_0$ is a
triangle with the vertices $A=0$, $B=1$ and $C=1/2+i(\sqrt{3}/2)$.
It is easily verified that $L(f_0)=1$. Suppose that $l_1$, $l_2$
and $l_3$ are the perpendiculars of sides $AB$, $BC$ and $AC$,
respectively. Then  $C_{f_0}=1/2+i\sqrt{3}/6$. Separate
$\mathbb{C}$ to three not disjoint partitions as the following:

$\Gamma_1$= The locations, which are on the right side or over
$l_1$ and top or over $l_3$.

$\Gamma_2$= The locations, which are on the left side or over
$l_1$ and top or over $l_2$.

$\Gamma_3$= The locations, which are under or over $l_2$ and under
or over $l_3$.

If we define $f(e)$ in $\Gamma_1$, then
$$\max\{|f(e)-A|,|f(e)-B|,|f(e)-C|\}=|f(e)-A|,$$
and so $$L(f)=\max\{1,2|f(e)-A|\}.$$ By the triangle inequality we
obtain $L(f)=2|f(e)-A|.$ But by Lemma \ref{right angle}, the
center of the triangle is at above of the base. Consequently by
Lemma \ref{3},
$$\min_{f(e)\in \Gamma_1}|f(e)-A|=|C_{f_0}-A|$$ and so
$$\min_{f(e)\in\Gamma_1}L(f)=2|C_{f_0}-A|=\frac{2}{\sqrt{3}}.$$
Similarly we obtain
$$\min_{f(e)\in \Gamma_2}L(f)=2|C_{f_0}-B|=\frac{2}{\sqrt{3}}$$ and
$$\min_{f(e)\in \Gamma_3}L(f)=2|C_{f_0}-C|=\frac{2}{\sqrt{3}}.$$
In all of the above cases, if $f$ is an extension of $f_0$, then
the following inequality is obtained
$$\frac{2}{\sqrt{3}}L(f_0)=\frac{2}{\sqrt{3}}\leq L(f).$$
Now suppose on the contrary that the inequality \eqref{2} is
satisfied for some $c<\dfrac{2}{\sqrt{3}}$. Thus we obtain
$$\frac{2}{\sqrt{3}}L(f_0)\leq L(f)\leq
cL(f_0)<\frac{2}{\sqrt{3}}L(f_0),$$ which is a contradiction.
These arguments imply that $\dfrac{2}{\sqrt{3}}$ is the best (the
least) possible constant, which satisfies the inequality
\eqref{2}.
\end{proof}

\begin{Note}
Note that in the assumption of Theorem \ref{4}, $l$ is assumed to
be greater than $\dfrac{1}{2}$. The reason is that if
$l<\dfrac{1}{2}$, then function $d$, defined in Theorem \ref{4} is
not a meter.
\end{Note}

The following result established a necessary and sufficient
condition for existing extension of Lipschitz, without increasing
Lipschitz number. It is obtained, by the proof of Theorem \ref{4}.

\begin{prop}\label{equal extension}
Suppose that $X_0$, $X$, $d$ and $l$ are those, as defined in
Theorem $\ref{4}$ and $f_0$ is a triangle on $X_0$. Then the
following statements are equivalent:
\begin{enumerate}
\item[(i)] $f_0$ has an angle that is greater than or equal to
$\sin^{-1}(\dfrac{1}{2l})$; \item[(ii)] $f_0$ has an extension $f$
on $X$  such that $L(f)=L(f_0)$.
\end{enumerate}
\end{prop}

\begin{proof}
$(i)\Rightarrow (ii)$. At first suppose that $f_0$ is normal and
$\theta$ is the largest angle of $f_0$. If
$\theta>\dfrac{\pi}{2}$, then define $f(e)=M_{f_0}$ and so the
result is obtained by the proof of Theorem \ref{4}. Otherwise if
$\theta\leq \dfrac{\pi}{2}$, then define $f(e)=C_{f_0}$. By Lemma
\ref{distance}, the distance between $C_{f_0}$ and the vertices,
is decreasing respect to $\theta$ and since $\theta\geq
\sin^{-1}(\dfrac{1}{2l})$, then $\dfrac{|C_{f_0}|}{l}\leq r$.
Consequently
$$
L(f)=\max\left\{L(f_0),\frac{|C_{f_0}|}{l}\right\}=r=L(f_0).
$$
So (ii) is obtained, whenever $f_0$ is normal. Now suppose that
$f_0$ is an arbitrary triangle on $X_0$. Then the conclusion is
obvious by the final argument in proof of first part of Theorem
\ref{4}.

$(ii)\Rightarrow (i)$. First suppose that $f_0$ is normal with the
vertices $A=0$, $B=r$ ($r>0$) and $C$. Suppose on the contrary
that all the angles are lower than $\sin^{-1}(\dfrac{1}{2l})$.
Consider three sets $\Gamma_1,\Gamma_2$ and $\Gamma_3$, defined as
in the proof of Theorem \ref{4}. Since $\theta\leq \dfrac{\pi}{2}$
then similar to proof of Theorem \ref{4}, for any extension $f$ of
$f_0$ on $X$ we have
\begin{equation}
\min_{f(e)\in \Gamma_1}L(f)=\min_{f(e)\in
\Gamma_2}L(f)=\min_{f(e)\in \Gamma_3}L(f)=\min_{f(e)\in
\mathbb{C}}L(f)=\frac{|C_{f_0}-A|}{l}=\frac{|C_{f_0}|}{l},
\end{equation}
and so
\begin{equation}\label{13}
L(f)\geq \frac{|C_{f_0}|}{l}.
\end{equation}
On the other since $\hat{C}<\sin^{-1}(\dfrac{1}{2l})$ and
$\hat{C}\leq \dfrac{\pi}{2}$, we obtain from Lemma \ref{distance}
that
\begin{equation}\label{12}
|C_{f_0}|>rl.
\end{equation}
Consequently \eqref{13} and \eqref{12} imply that
$$
L(f)\geq \frac{|C_{f_0}|}{l}>r=L(f_0).
$$
It follows that $f_0$ has not any extension $f$ with
$L(f)=L(f_0)$. Now suppose that $f_0$ is an arbitrary triangle
such that all of its angles are lower than
$\sin^{-1}(\dfrac{1}{2l})$, and $f$ is an arbitrary extension of
$f_0$ to $X$. Consider the transformations, done on $f_0$ to
obtain ${f_0}^{\sharp}$. Take $M\in \mathbb{C}$ to be the point,
obtained by doing these transformations on $f(e)$. Note that all
the angles of ${f_0}^{\sharp}$ are also lower than
$\sin^{-1}(\dfrac{1}{2l})$. Then from the first part of the proof,
${f_0}^{\sharp}$ has not any extension with the same Lipschitz
number of ${f_0}^{\sharp}$. Thus if $f^{\sharp}$ is an extension
of ${f_0}^{\sharp}$ such that $f^{\sharp}(e):=M$, then
$L({f_0}^{\sharp})<L(f^{\sharp})$. Now Proposition \ref{lem}
implies that
$$
L(f_0)=L({f_0}^{\sharp})<L(f^{\sharp})=L(f).
$$
Thus the proof is completed.
\end{proof}

The following result is immediately obtained from Proposition
\ref{equal extension}.

\begin{cor}
Suppose that $X_0$, $X$ and $d$ are those, as defined in Theorem
$\ref{4}$, $l\geq\dfrac{1}{\sqrt{3}}$ and $f_0$ is a triangle on
$X_0$. Then $f_0$ has an extension $f$ on $X$  such that
$L(f)=L(f_0)$.
\end{cor}

In the next result, we extend Proposition \ref{equal extension}
for some special tetragon.

\begin{prop}\label{tetragone}
Suppose that $X_0$, $X$, $d$ and $l$ are those, as defined in
Theorem $\ref{4}$ and $f_0$ is a tetragon on $X_0$. Moreover
suppose that $f_0$ has two angles which are greater than or equal
to $\sin^{-1}(\dfrac{1}{2l})$ and are also against of each other.
Then $f_0$ has an extension  $f$ on $X$ such that $L(f)=L(f_0)$.
\end{prop}

\begin{proof}
Suppose that the range of $f_0$ is tetragon $ABCD$ and two angles
which are  greater than or equal to $\sin^{-1}(\dfrac{1}{2l})$ and
against of each other are $A$ and $C$. We denote by $g_0$ and
$h_0$, the triangles $ABD$ and $BCD$, respectively. In fact $g_0$
and $h_0$ are obtained via the eliminating $f_0^{-1}(C)$ and
$f_0^{-1}(A)$, from the domain of $f_0$. Define $f(e)$ to be the
middle point of the diagonal $BD$. We prove that such a this $f$
is the desired extension for $f_0$. To that end, suppose that $g$
and $h$ are the extensions of $g_0$ and $h_0$ on $X$ respectively,
defined as $g(e)=h(e):=f(e)$. By Proposition \ref{equal extension}
we have $L(g)=L(g_0)$ and $L(h)=L(h_0).$ Consequently
$$
\max\{\dfrac{1}{l}|f(e)-A|,\dfrac{1}{l}|f(e)-B|,\dfrac{1}{l}|f(e)-D|\}\leq
L(g_0)
$$
and
$$
\max\{\dfrac{1}{l}|f(e)-B|,\dfrac{1}{l}|f(e)-C|,\dfrac{1}{l}|f(e)-D|\}\leq
L(h_0).
$$
But we have
$$L(g_0)\leq L(f_0)\;\;\;and\;\;\; L(h_0)\leq L(f_0).$$
Thus
\begin{eqnarray*}
L(f)&=&\max\{L(f_0),\dfrac{1}{l}|f(e)-A|,\dfrac{1}{l}|f(e)-B|,\dfrac{1}{l}|f(e)-C|,\dfrac{1}{l}|f(e)-D|\}\\&\leq&
\max\{L(f_0),L(g_0),L(h_0)\}\\&=& L(f_0).
\end{eqnarray*}
Therefore $L(f)=L(f_0)$.
\end{proof}

\begin{cor}
Suppose that $X_0$, $X$ and $d$ are those, as defined in Theorem
$\ref{4}$, $l=\dfrac{1}{2}$ and $f_0$ is a tetragon on $X_0$. If
$f_0$ has two  wide (or right) angles which are against of each
other, then $f_0$ has an extension $f$ on $X$ such that
$L(f)=L(f_0)$.
\end{cor}

\begin{cor}
Suppose that $X_0$, $X$ and $d$ are those, as defined in Theorem
$\ref{4}$, $l=\dfrac{1}{2}$. Then every parallelogram on $X_0$ has
an extension on $X$, without increasing Lipschitz number.
\end{cor}

In the next result, we generalize Theorem \ref{4} for some special
regular $n-$polygons.

\begin{thrm}\label{10}
Suppose that $X_0$, $X$, $d$ and $l$ are those, as defined in
Theorem $\ref{4}$. Moreover suppose that $f_0$ is a regular
$n$-polygon on $X_0$ with $n\geq 3$. Then $f_0$ has an extension
$f$ on $X$ such that
$$
L(f)\leq \frac{2}{\sqrt{3}}L(f_0).
$$
Moreover the best (least) possible constant for the above
inequality is $\dfrac{2}{\sqrt{3}}$.
\end{thrm}

\begin{proof}
It is easily verified that $L(f_0)$ is exactly the length of the
largest diameter of $f_0$. We denote it by $D$. Note that for the
case where $n=3$, $D$ is the length of the largest side. Define
$f(e):=C_{f_0}$. Thus
\begin{equation}\label{16}
L(f)=\max\{D,\frac{E}{l}\}.
\end{equation}
Since $l\geq\dfrac{1}{2}$, thus  $\dfrac{E}{l}\leq 2E$ and so
$L(f)=\max\{D,2E\}.$ But from the triangle inequality we obtain
$2E\geq D$. Consequently
\begin{equation}\label{15}
\frac{L(f)}{L(f_0)}\leq \frac{2E}{D}.
\end{equation}
Denote by $a$ the length of sides of $f_0$. By some calculations
we obtain
$$
2E^2-2E^2\cos(\frac{2\pi}{n})=a^2.
$$
Thus
$$ E=\frac{a}{\sqrt{2(1-\cos(\frac{2\pi}{n}))}}=\frac{a}{2\sin(\frac{\pi}{n})}.
$$
Now consider the following two cases:
\begin{enumerate}
\item[(i)] Suppose that $n$ is odd. Take $A$ to be one of the
arbitrary vertices of $f_0$. Then there are two vertices of $f_0$,
denoted by $B$ and $C$ such that their distances until $A$ is
maximal. In triangle $ABC$, we have $A=\dfrac{\pi}{n}$,
$|AB|=|AC|=D$ and $|BC|=a$. Thus
$$
2D^2-2D^2\cos(\frac{\pi}{n})=a^2
$$
and so
$$
D=\frac{a}{2\sin(\frac{\pi}{2n})}.
$$
Consequently
\begin{eqnarray*}
\frac{2E}{D}=\frac{2\sin(\frac{\pi}{2n})}
{\sin(\frac{\pi}{n})}=\frac{1}{\cos(\frac{\pi}{2n})}
\end{eqnarray*}
\item[(ii)] Suppose that $n$ is even. It is easily verified that
the largest diameter of $f_0$ crosses $C_{f_0}$. It follows that
$$D=2E=\frac{a}{\sin(\frac{\pi}{n})}$$ and consequently
$\dfrac{2E}{D}=1$.
\end{enumerate}
Thus \eqref{15} implies that
$$
\frac{L(f)}{L(f_0)}\leq \max\left\{1,\max_{n\geq 3}\{\frac{1}{\cos(\frac{\pi}{2n})}\}\right\}\leq \frac{2}{\sqrt{3}}.
$$
Moreover by Theorem \ref{4}, the best (least) possible constant
for this inequality is $\dfrac{2}{\sqrt{3}}$.
\end{proof}

\begin{cor}\label{equal n}
Suppose that $X_0$, $X$, $d$ and $l$ are those, as defined in
Theorem $\ref{10}$ and $f_0$ is a regular $n$-polygon on $X_0$. If
$n$ is even or $n\geq\dfrac{\pi}{2\cos^{-1}(\frac{1}{2l})}$, then
$f_0$ has an extension to $X$, without increasing Lipschitz
number.
\end{cor}

\begin{proof}
Define $f(e):=C_{f_0}$ and suppose that $D$ and $E$ are as in
Theorem \ref{10}. If $n$ is even, then by Theorem \ref{10} we have
$L(f)=L(f_0)$. Now suppose that
$n\geq\dfrac{\pi}{2\cos^{-1}(\dfrac{1}{2l})}$. By some easy
calculations we obtain
\begin{eqnarray*}
 \dfrac{E/l}{D}= \left\lbrace
 \begin{array}{c l}
  \dfrac{1}{2l\cos(\dfrac{\pi}{n})}& \hbox{\;\;  if $n$ is odd}\\
  \dfrac{1}{2l}& \hbox{\;\;  if $n$ is even}

 \end{array}\right. \leq 1
\end{eqnarray*}
Consequently (\ref{16}) implies that $L(f)=D=L(f_0)$.
\end{proof}

\section{\bf Tietze-Lipschitz extension}

In topology, the Tietze extension theorem (also known as the
Tietze-Urysohn-Brouwer extension theorem) states that every
continuous real-valued function on a closed subset of a normal
topological space can be extended to the entire space, preserving
supremum norm. Also in \cite[Theorem 1.5.6]{W}, it has been shown
that if $(X,d)$ is a metric space and $X_0\subseteq X$, then every
real-valued function $f\in Lip_1 X_0$ has an extension $f\in
Lip_1X$ without increasing Lipschitz number, in particular
$\|f_0\|_1=\|f\|_1.$ Moreover, it has been proved that every
complex-valued function $f\in Lip_1 X_0$ has an extension $f\in
Lip_1X$ such that $L(f)\leq\sqrt{2}\;L(f_0)$ and
$\|f_0\|_\infty=\|f\|_\infty$. In fact for the complex version, it
has not been constructed any extended Lipschitz function,
preserving Lipschitz norm.

Here, we try to generalize Tietze extension theorem for
complex-valued Lipschitz functions. In fact we create conditions
under which every $f\in Lip_\alpha X_0$, where $X_0=X-\{e\}$ and
$\alpha>0$, can be extended to a function $f\in Lip_\alpha X$,
preserving Lipschitz norm.

\begin{thrm}\label{d(y,e)>d(x,y)}
Let $(X,d)$ be a metric space, $e\in X$ and $X_0=X-\{e\}$ and
$\alpha>0$. Moreover suppose that for all $x,y\in X_0$ we have
\begin{equation}\label{8}
d(y,e)\geq d(x,y).
\end{equation}
Then for each $f_0\in Lip_{\alpha}X_0$, there exists an extension
$f\in Lip_{\alpha}X$ such that $p_{\alpha}(f)=p_{\alpha}(f_0)$ and
$\|f\|_{\infty}=\|f_0\|_{\infty}$. In particular $\|
f\|_{\alpha}=\|f_0\|_{\alpha}$.
\end{thrm}

\begin{proof}
We should only define $f(e)$. For each finite subset $F$ of $X$,
let
$$
M_F(f):=\frac{\sum_{x\in F}f(x)}{|F|},
$$
where $|F|$ is the number of elements of $F$. It is easily
verified that if $f$ is bounded, then
$$
|M_F(f)|\leq\|f\|_\infty.
$$
Consider the net $\{M_F(f_0)\}_{F\in\mathcal{F}}$, where
$\mathcal{F}$ is the collection of finite subsets of $X_0$,
directed by upward inclusion. Since for every $F\in{\mathcal F}$
\begin{equation}\label{11}
|M_F(f_0)|\leq\|f_0\|_{\infty},
\end{equation}
the net $\{M_F(f_0)\}_{F\in\mathcal{F}}$ is in the compact set
$\{z\in \mathbb{C} :\;\;|z|\leq \|f_0\|_{\infty}\}$. Thus it
contains a convergent subnet, called
$\{M_F(f_0)\}_{F\in\mathcal{F}}$ again by abuse of notation.
Define
$$f(e)=  \lim_{F}M_F(f_0).$$
For every finite subset $F\subseteq X_0$ and $y\in X_0$ we have
\begin{eqnarray*}
|M_F(f_0)-f_0(y)| &\leq& \frac{1}{|F|}\sum_{x\in
F}|f_0(x)-f_0(y)|\\&\leq&\frac{1}{|F|}\sum_{x\in
F}p_{\alpha}(f_0)d^{\alpha}(x,y).
\end{eqnarray*}
By the inequality \eqref{8}, for such $F$ and $y$ we obtain
$$
\frac{|M_F(f_0)-f_0(y)|}{d^{\alpha}(e,y)}\leq\frac{1}{|F|}\sum_{x\in
F}p_{\alpha}(f_0)(\frac{d(x,y)}{d(e,y)})^{\alpha}\leq
p_{\alpha}(f_0).
$$
Thus
$$ \sup_{y\in
M}\frac{|f(e)-f_0(y)|}{d^{\alpha}(e,y)}=  \sup_{y\in
M}\lim_{F}\frac{|M_{F}(f_0)-f_0(y)|}{d^{\alpha}(e,y)}\leq
p_{\alpha}(f_0)
$$
and consequently
$$p_{\alpha}(f)=\max\left\{p_{\alpha}(f_0),\sup_{y\in
M}\{\frac{|f(e)-f_0(y)|}{d^{\alpha}(e,y)}\}\right\}=p_{\alpha}(f_0).
$$
That $\|f\|_{\infty}=\|f_0\|_{\infty}$ is obvious. Therefore
$$\|f\|_{\alpha}=\|f_0\|_{\alpha}$$
and so the result is proved.
\end{proof}

Note that inequality \eqref{8} is a sufficient condition for the
existence of an extension with preserving Lipschitz norm. But
Proposition \ref{equal extension} confirms that it is not a
necessary condition.

Recall that if $(X,d)$ is a metric space and $M\subseteq X$, then
the diameter of $M$ is defined as
$$diam(M):=\sup_{x,y\in M}d(x,y).$$ Moreover for each $x\in X$
$$
d(M,x)=\inf_{y\in M}d(x,y).
$$
Now the following corollary is obtained from Theorem
\ref{d(y,e)>d(x,y)}, immediately.

\begin{cor}
Let $(X,d)$ be a metric space, $e\in X$ and $X_0=X-\{e\}$ and
$\alpha>0$. Also suppose that
\begin{equation}\label{8.1}
d(X_0,e)\geq diam(X_0).
\end{equation}
Then for each function $f_0\in Lip_{\alpha}(X_0)$, there exists an
extension $f\in Lip_{\alpha}(X)$ such that
$p_{\alpha}(f)=p_{\alpha}(f_0)$ and $\|
f\|_{\infty}=\|f_0\|_{\infty}$. In particular $f\in
Lip_{\alpha}(X)$ and $\|f\|_\alpha=\|f_0\|_\alpha$.
\end{cor}

\begin{remark}\rm
Note that inequality \eqref{8} is a necessary but not sufficient
condition for inequality \eqref{8.1}. For example let
$X=\{x,y,z,e\}$ and $M=\{x,y,z\}$. Define the metric $d$ on $X$ by
$d(x,e)=d(z,e)=d(x,z)=5$, $d(y,e)=4$ and $d(x,y)=d(y,z)=3$. Then
$(X,d)$ is a metric space, satisfying \eqref{8}. But
$$4=d(M,e)\ngeq diam(M)=5.$$
\end{remark}

\begin{prop}\label{d(y,z)>d(x,y)}
Let $(X,d)$ be a metric space, $X_0\subseteq X$ and $\alpha>0$.
Moreover suppose that the meter $d$ satisfies the following
inequality
$$d(x,y)\geq d(y,z),$$
for all $x\in X-X_0$ and $y,z\in X_0$. Then every $f_0\in
Lip_\alpha(X_0)$ has an extension $f$ on $X$ such that
$p_{\alpha}(f)=p_{\alpha}(f_0)$ and
$\|f\|_{\infty}=\|f_0\|_{\infty}$. In particular $f\in
Lip_{\alpha}(X)$ and $\|f\|_\alpha=\|f_0\|_\alpha$.
\end{prop}

\begin{proof}
Define $f$ such as that defined in Theorem \ref{d(y,e)>d(x,y)}.
Then by some arguments similar to the proof of Theorem
\ref{d(y,e)>d(x,y)} the conclusion is obtained.
\end{proof}

\begin{cor}
Let $(X,d)$ be a metric space, $X_0\subseteq X$ and $\alpha>0$.
Moreover suppose that $d$ is a metric on $X$ such that
$$d(X-X_0,X_0)\geq diam(X_0).$$
Then every $f_0\in Lip_\alpha(X_0)$ has an extension $f$ on $X$
such that $p_{\alpha}(f)=p_{\alpha}(f_0)$ and
$\|f\|_{\infty}=\|f_0\|_{\infty}.$ In particular $f\in
Lip_{\alpha}(X)$ and $\|f\|_\alpha=\|f_0\|_\alpha$.
\end{cor}

Our next results rely on Helly`s theorem \cite{H}, which states
that if every triple of sets in a collection of closed convex
subsets of $\mathbb{C}$ have nonempty intersection, then the whole
collection has nonempty intersection. First, we require the
following elementary lemma. The proof is straightforward and so is
left to the reader.

\begin{lem}\label{circle}
Let $P_1=x_1+iy_1$ and $P_2=x_2+iy_2$ be two points in
$\mathbb{C}$, $K>0$ and
$$I:=\{z\in \mathbb{C}:|\frac{z-P_1}{z-P_2}|=K\}.$$
\begin{enumerate}
\item[(i)] If $K=1$, then $I$ is the perpendicular of the segment
$P_1P_2$. \item[(ii)] If $K\neq 1$, then $I$ is a circle with
center $\alpha+i\beta$ and radios $R$, such that
$$
\alpha=\frac{K^2x_1-x_2}{K^2-1}\;\;\;,\;\;\;\beta=\frac{K^2y_1-y_2}{K^2-1}
\;\;\;\hbox{and}\;\;\; R=\frac{K}{|K^2-1|}|P_1-P_2|
$$
\end{enumerate}
\end{lem}

\begin{prop}\label{14}
Let $(X,d)$ be a metric space, $ X_0=\{x_1,x_2,...,x_n\}\subseteq
X$, $e\in X-X_0$ and $\alpha>0$. Moreover for each $j=1,2,...,n-1$
let
$$D_j=\left\{z\in \mathbb{C}:|\frac{z-f(x_j)}{z-f(x_{j+1})}|\leq \frac{d_j}{d_{j+1}}\right\},$$
and $$D_n=\{z\in \mathbb{C}:|z-f(x_n)|\leq d_n p_\alpha(f)\},$$
where $d_j:=d(e,x_j)$.
\begin{enumerate}
\item[(i)] If for every three distinct elements $m$, $k$ and $l$
of the set $\{1,2,...,n-1\}$ we have $D_k\cap D_l\cap D_m\cap
D_n\neq \emptyset$ then each function $f_0\in Lip_\alpha(X_0)$ has
an extension $f\in Lip_{\alpha}(X_0\cup \{e\})$ such that
$p_\alpha(f)=p_\alpha(f_0)$. \item[(ii)] If there exists $1\leq
j\leq n$ such that $|f(e)-f_0(x_j)|>d_jp_\alpha(f_0)$ then
$$p_\alpha(f)>p_\alpha(f_0).$$
\end{enumerate}
\end{prop}

\begin{proof}
(i). Set $B_j=D_j\cap D_n$, for every $j=1,2,...,n-1$. By Lemma
\ref{circle} every $D_j$ is either a disc or a Half plane in
$\mathbb{C}$. Thus each $B_j$ is closed and convex in
$\mathbb{C}$. By the assumption and Helly`s theorem we have
$$\cap_{j=1}^{n-1}B_j=\cap_{j=1}^nD_j\neq \emptyset.$$ Now define
$f(e)=M$, such that $M$ is an arbitrary point in
$\cap_{j=1}^{n}D_{j}$. Then we have
$$\frac{|f(e)-f_0(p_1)|}{d_1^{\alpha}}\leq \frac{|f(e)-f_0(p_2)|}
{d_2^{\alpha}}\leq ...\leq
\frac{|f(e)-f_0(p_{n})|}{d_{n}^{\alpha}}\leq p_{\alpha}(f_0).$$
This implies that
$$p_{\alpha}(f)=\left\{\max_{j\in \{1,2,...,n\}}\frac{|f(e)-f_0(p_j)|}
{d_j^{\alpha}},p_{\alpha}(f_0)\right\}=p_{\alpha}(f_0).$$ (ii).
This is obvious.
\end{proof}

Let $(X,d)$ be a metric space, $X-X_0=\{e\}$, $\alpha>0$ and
$f_0\in Lip_{\alpha}X_0$. For $x\in X$ and $r>0$, let
$$B_{x,f_0}=\{z\in\Bbb C:\;\; |z-f_0(x)|\leq d(x,e)p_\alpha(f_0)\},$$
$$B_{0,f_0}=\{z\in\Bbb C:\;\; |z|\leq\|f_0\|_\infty\}.$$
Note that all these spaces are closed and convex in $\Bbb C$.

\begin{prop}
Let $(X,d)$ be a metric space, $X-X_0=\{e\}$, $\alpha>0$ and
$f_0\in Lip_{\alpha}X_0$. Then $f_0$ has an extension $f$ to $X$
such that
\begin{enumerate}
\item[(i)] $p_{\alpha}(f)=p_{\alpha}(f_0)$ if and only if for
every $x,y,z\in X_0$
$$B_{x,f_0}\cap B_{y,f_0} \cap B_{z,f_0}\neq \emptyset.$$
\item[(ii)] $\|f\|_{\alpha}=\|f_0\|_{\alpha}$ if and only if for
every $x,y,z\in X_0$
$$B_{0,f_0}\cap B_{x,f_0}\cap B_{y,f_0}\cap B_{z,f_0}\neq \emptyset.$$
\end{enumerate}
\end{prop}

\begin{proof}
(i). Suppose that $B_{x,f_0}\cap B_{y,f_0} \cap B_{z,f_0}\neq
\emptyset$, for all $x,y,z\in X_0$. By Helly`s theorem
$$
\cap_{x\in X_0}B_{x,f_0}\neq \emptyset
$$
Now choose $f(e)\in \cap_{x\in X_0}B_{x,f_0}$. Thus
$$\frac{|f(e)-f_0(x)|}{d(x,e)}\leq p_{\alpha}(f_0)$$ for every $x\in
X_0$. It follows that $p_{\alpha}(f)=p_{\alpha}(f_0)$.

(ii). Let $D_{x,f_0}=B_{x,f_0}\cap B_{0,f_0}$ for each $x\in X_0$
and use again Helly`s theorem.
\end{proof}

\section{\bf Urysohn-Lipschitz extension}

Two disjoint subsets $A$ and $B$ of the metric space $(X,d)$, are
said to be separated by a function if there exists a continuous
function $f:X\rightarrow [0,1]$ such that $f(A)=0$ and $f(B)=1$.
Any such function is called a Urysohn function for $A$ and $B$.
Urysohn's lemma states that a topological space is normal if and
only if any two disjoint closed subsets can be separated by a
continuous function. In this section, we generalize Urysohn's
lemma as follows; we present a necessary and sufficient condition
for that Lipschitz functions separate two subsets $A$ and $B$.

We commence with following proposition.

\begin{prop}\label{exa}
Let $(X,d)$ be a metric space. Suppose that $A$ and $B$ are two
subsets of $X$ such that $d(A,B)>0$. Then there is a function
$f\in Lip_{1}X$, satisfying the following conditions:
\begin{enumerate}
\item[(i)] $0\leq f(x)\leq 1$\;\;\;\;\;\;\;\;$(x\in X)$,
\item[(ii)] $f(A)=0\;$ and $\;f(B)=1$, \item[(iii)] $L(f)\leq
\dfrac{1}{d(A,B)}$.
\end{enumerate}
\end{prop}

\begin{proof}
Define the function $f$ on $X$ as
$$
f(x)=\frac{d(x,A)}{d(x,A)+d(x,B)}.
$$
First note that since $d(A,B)>0$, then $A\cap B=\emptyset$ and so
$f$ is well defined. Also one can easily verified that $f$
satisfies $(i)$ and $(ii)$. Thus we only prove part $(iii)$. For
all $x,y\in X$ with $x\neq y$ we have
\begin{eqnarray*}
|f(x)-f(y)|&=&\left|\frac{d(x,A)}{d(x,A)+d(x,B)}-\frac{d(y,A)}{d(y,A)+d(y,B)}\right|\\&=&
\frac{|d(x,A)d(y,B)-d(y,A)d(x,B)|}{[d(x,A)+d(x,B)][d(y,A)+d(y,B)]}\\&\leq&
\frac{d(x,A)|d(y,B)-d(x,B)|+d(x,B)|d(y,A)-d(x,A)|}{[d(x,A)+d(x,B)][d(y,A)+d(y,B)]}\\&\leq&
\frac{d(x,A)d(x,y)+d(x,B)d(x,y)}{[d(x,A)+d(x,B)][d(y,A)+d(y,B)]}\\&\leq&
\frac{d(x,y)}{d(y,A)+d(y,B)}.
\end{eqnarray*}
Since
$$
d(y,A)+d(y,B)\geq d(A,B),
$$
we obtain
$$
|f(x)-f(y)|\leq \frac{1}{d(A,B)}d(x,y).
$$
It follows that $L(f)\leq\dfrac{1}{d(A,B)}$ and so $(iii)$ is
obtained.
\end{proof}

If $d$ is a metric on $X$, then $d^\alpha$ is also a metric on
$X$, whenever $0<\alpha<1$. Thus the following result is obtained
from Proposition \ref{exa}.

\begin{cor}\label{separator}
Let $(X,d)$ be a metric space and $A, B$ are two subsets of $X$
such that $d(A,B)>0$. Then for each $0<\alpha\leq 1$, there is a
function $f\in Lip_{\alpha}X$, satisfying the following
conditions:
\begin{enumerate}
\item[(i)] $0\leq f\leq 1$, \item[(ii)] $f(A)=0$ and $f(B)=1$,
\item[(iii)] $L(f)\leq \dfrac{1}{d^{\alpha}(A,B)}$.
\end{enumerate}
\end{cor}

\begin{proof}
By the proof of Proposition \ref{exa}, function $f$ defined as
$$
f(x)=\frac{d^\alpha(x,A)}{d^\alpha(x,A)+d^\alpha(x,B)}\;\;\;\;\;\;\;\;\;\;\;\;\;(x\in
X)
$$
satisfies all the conditions.
\end{proof}

\begin{remark}\label{r1}\rm
Let $(X,d)$ be a metric space and $A$ and $B$ be two disjoint
subsets of $X$ such that $A$ is compact and $B$ is closed. We show
that $d(A,B)>0$. To that end, suppose on the contrary that
$d(A,B)=0$. Then there exists sequence $(x_n)_n$ in $A$ such that
$$0\leq d(x_n,B)<1/n\;\;\;\;\;\;\;(n\in\Bbb N).$$ Since $A$ is
compact, there is subsequence $x_{n_k}$ of $(x_n)_n$, converging
to $x$, for some $x\in A$. It follows that
$$
d(x,B)=\lim_{k}d(x_{n_k},B)=0,
$$
which implies that $x\in \overline{B}=B$. Thus $x\in A\cap B$,
that is a contradiction. It follows that $d(A,B)>0$. Therefore
Proposition \ref{exa} and Corollary \ref{separator} hold, whenever
$A$ is compact and $B$ is closed.
\end{remark}

Let $X$ be a nonempty set and $\alpha>0$. Suppose that $A$ and $B$
are two disjoint subsets of $X$. We say that $Lip_\alpha X$
separates $A$ and $B$ if there exists $f\in Lip_\alpha X$ and
different scalers $\lambda_1$ and $\lambda_2$ such that
$f(A)=\lambda_1$ and $f(B)=\lambda_2$.

Here we state the main result of this section.

\begin{thrm}\label{separating condition}(Lipschitz-Urysohn)
Let $(X,d)$ be a metric space, $A,B\subseteq X$ and $0<\alpha\leq
1$. Then the following assertions are equivalent:
\begin{enumerate}
\item[(i)] $d(A,B)>0$, \item[(ii)] $Lip_{\alpha}(X)$ separates $A$
and $B$.
\end{enumerate}
\end{thrm}

\begin{proof}
$(i)\Rightarrow (ii)$. It is obtained from Corollary
\ref{separator}.

$(ii)\Rightarrow (i)$. Suppose that there exists $f\in Lip_\alpha
X$ and different scalers $\lambda_1$ and $\lambda_2$ such that
$f(A)=\lambda_1$ and $f(B)=\lambda_2$. Thus there exists $M>0$
such that
\begin{equation}\label{e2}
|f(x)-f(y)|\leq M d^{\alpha}(x,y),
\end{equation}
for all $x,y\in X$ with $x\neq y$. Suppose on the contrary that
$d(A,B)=0$. It follows that there are sequences $(x_n)$ and
$(y_n)$ in $A$ and $B$, respectively such that $$d(A,B)=\lim_{n\to
\infty}d(x_n,y_n)=0.$$ Now inequality \eqref{e2} implies that
$$0<|f(x_n)-f(y_n)|\leq M d^{\alpha}(x_n,y_n)\to 0,$$ which is a contradiction.
Consequently $d(A,B)>0$ and so $(i)$ is obtained.
\end{proof}

\begin{remark}\rm
\begin{enumerate}
\item We notify that in the proof of implication $(ii)\Rightarrow
(i)$ of Theorem \ref{separating condition}, we did not use the
condition $0<\alpha\leq 1$. Thus $(ii)\Rightarrow (i)$ is valid
for any $\alpha>0$. \item Note that, Theorem \ref{separating
condition} is not necessarily true when $\alpha>1$. For example
let $X=\mathbb{R}$ with the euclidian metric and $\alpha=2$. It is
not hard to see that $Lip_{2}(\mathbb{R})=Cons(\mathbb{R})$. It
follows that $Lip_{2}(\Bbb R)$ does not separate any two disjoint
subsets of $\mathbb{R}$. \item Let $(X,d)$ be a metric space and
$0<\alpha\leq 1$. Remark \ref{r1} implies that $Lip_\alpha X$
separates two disjoint subsets $A$ and $B$ such that one of them
is compact and another one is closed.
\end{enumerate}
\end{remark}

\begin{Note}
Let $(X,d)$ be a metric space and $0<\alpha\leq 1$. In \cite[page
4]{W}, it has been claimed that for all disjoint closed subsets
$A$ and $B$ of $X$, there is a Lipschitz function
$f:X\rightarrow\Bbb R$ which is constantly zero on $A$ and
constantly 1 on $B$. In fact it was mentioned that the existence
of this function is a kind of metric version of Urysohn's Lemma.
Whereas in Theorem \ref{separating condition}, we showed that the
existence of this function is equivalent to $d(A,B)>0$. It is
understood that this false claim may obtain from this wrong
thinking that $d(A,B)>0$, for all closed and disjoint subsets $A$
and $B$ of $X$. We provide an example to demonstrate that this is
not the case in general. Take $X$ to be $\mathbb{R}$, endowed with
the Euclidean metric, $0<\alpha\leq 1$ and suppose that
$$A=\{n:n\in
\mathbb{N}\}\;\; and \;\; B=\{n+\frac{1}{n}:n\in \mathbb{N}\}.$$
Then $A$ and $B$ are disjoint and closed subsets of $\mathbb{R}$
such that $d(A,B)=0$. In fact $$d(A,B)\leq d(x_n,y_n)\to 0.$$ Now
Theorem \ref{separating condition} implies that $Lip_{\alpha}X$
does not separate such a $A$ and $B$.
\end{Note}

\end{document}